\documentclass[12pt,reqno]{amsart}
\usepackage{amsmath}
\usepackage{amssymb}
\usepackage[left=3cm,top=3cm,right=3cm,bottom=3cm]{geometry}
\usepackage{epsfig}

\begin{document}
\newtheorem{theorem}{Theorem}[section]
\newtheorem{lemma}[theorem]{Lemma}
\newtheorem{definition}[theorem]{Definition}
\newtheorem{conjecture}[theorem]{Conjecture}
\newtheorem{proposition}[theorem]{Proposition}
\newtheorem{algorithm}[theorem]{Algorithm}
\newtheorem{corollary}[theorem]{Corollary}
\newtheorem{observation}[theorem]{Observation}
\newtheorem{problem}[theorem]{Open Problem}
\newcommand{\noin}{\noindent}
\newcommand{\ind}{\indent}
\newcommand{\om}{\omega}
\newcommand{\pp}{\mathcal P}
\newcommand{\ppp}{\mathfrak P}
\newcommand{\N}{{\mathbb N}}
\newcommand{\LL}{\mathbb{L}}
\newcommand{\R}{{\mathbb R}}
\newcommand{\E}{\mathbb E}
\newcommand{\Prob}{\mathbb{P}}
\newcommand{\eps}{\varepsilon}

\newcommand{\Ss}{{\mathcal S}}
\newcommand{\Nn}{{\mathcal N}}

\title{Revolutionaries and spies on random graphs}

\author{Dieter Mitsche}
\address{Department of Mathematics, Ryerson University, Toronto, ON, Canada}
\email{\tt dmitsche@ryerson.ca}

\author{Pawe\l{} Pra\l{}at}
\address{Department of Mathematics, Ryerson University, Toronto, ON, Canada}
\email{\tt pralat@ryerson.ca}

\thanks{The authors gratefully acknowledge support from NSERC and MPrime}

\maketitle

\begin{abstract}
Pursuit-evasion games, such as the game of Revolutionaries and Spies, are a simplified model for network security. In the game we consider in this paper, a team of $r$ revolutionaries tries to hold an unguarded meeting consisting of $m$ revolutionaries. A team of $s$ spies wants to prevent this forever. For given $r$ and $m$, the minimum number of spies required to win on a graph $G$ is the spy number $\sigma(G,r,m)$. We present asymptotic results for the game played on random graphs $G(n,p)$ for a large range of $p = p(n), r=r(n)$, and $m=m(n)$. The behaviour of the spy number is analyzed completely for dense graphs (that is, graphs with average degree at least $n^{1/2+\eps}$ for some $\eps > 0$). For sparser graphs, some bounds are provided. 
\end{abstract}

\section{Introduction}\label{sec:intro}

In the past few years, several problems from applications related to the structure of modern networks in the real-world emerged. In these problems the behaviour of the agents can be probabilistic, decentralized and even selfish or antagonistic. This is one of the reasons why the field of graph searching is nowadays rapidly expanding. Several new models, problems or approaches have appeared relating it to diverse fields such as random walks, game theory, logic, probabilistic analysis, complex networks, motion planning, and distributed computing. Surprising and not yet widely circulated results have been found in the last few years that had consequences for the whole field. For more details see, for example,~\cite{al,fomin}.

Suppose that unsupervised intruders are on the vertices of a network, and suppose that they travel between adjacent vertices. The intruders could represent viruses, hackers, or some other malicious agents that want to attack some vertex in the network. A set of searchers is attempting to prevent them from doing this. Although placing a searcher on each vertex protects the whole network, it is desired to find the minimum number of searchers required to do it. A motivation for minimizing the number of searchers comes from the fact that fewer searchers require fewer resources. Networks that require a smaller number of searchers may be viewed as more secure than those where many searchers are needed.

A pursuit-evasion game may be viewed as a simplified model for such network security problems. In this paper, we study the game of Revolutionaries and Spies invented by J\'{o}zsef Beck in the mid-1990s (according to~\cite{CSW}) that is played on a fixed graph $G$. There are two players: a team of $r$ \emph{revolutionaries} and a team of $s$ \emph{spies}. The revolutionaries want to arrange a one-time meeting consisting of $m$ revolutionaries free of oversight by spies; the spies want to prevent this from happening. The revolutionaries start by occupying some vertices as their initial positions; more than one revolutionary is allowed to occupy some vertex. After that the spies do the same; they can start from some vertices already occupied by revolutionaries or choose brand new vertices. In each subsequent round, each revolutionary may move to an adjacent vertex or choose to stay where he is, and then each spy has the same option. This is a perfect information game, that is, both players see each other's locations and know all possible moves. Moreover, we assume that the players are perfect, that is, they can analyze the game completely and play optimally. For more on combinatorial games see, for example,~\cite{anw}.

A \emph{meeting} is a set of at least $m$ revolutionaries occupying one vertex; a meeting is \emph{unguarded} if there is no spy at that vertex. The revolutionaries win if at the end of some round there is an unguarded meeting. On the other hand, the spies win if they have a strategy to prevent this forever. For given $r$ and $m$, the minimum number of spies required to win on a graph $G$ is the \emph{spy number} $\sigma(G,r,m)$. Since $\min\{ |V(G)|, \lfloor r/m \rfloor \}$ meetings can be initially formed, at least that many spies are needed to have a chance to win. On the other hand, $r-m+1$ spies can create a matching with $r-m+1$ distinct revolutionaries and follow them during the whole game, preventing any unguarded meeting from taking place. If $|V(G)|<r-m+1$, this can be clearly improved since occupying all vertices clearly does the job as well. We thus get the following trivial bounds on the spy number:
$$
\min\{ |V(G)|, \lfloor r/m \rfloor \} \le \sigma(G,r,m) \le \min \{ |V(G)|, r-m+1\}.
$$
It is known that the lower bound is sufficient when $G$ is a tree, and at most one additional spy is needed to win on any unicyclic graph~\cite{CSW}. On the other hand, the upper bound can be obtained for some chordal and bipartite graphs (such as hypercubes, for example)---see~\cite{BCPWZ} for more. Moreover, grid-like graphs were studied in~\cite{HS}. 

\bigskip

In this paper, we consider the spy number in binomial random graphs. The \emph{random graph} $G(n,p)$ consists of the probability space $(\Omega, \mathcal{F}, \Prob)$, where $\Omega$ is the set of all graphs with vertex set $[n]=\{1,2,\dots,n\}$, $\mathcal{F}$ is the family of all subsets of $\Omega$, and for every $G \in \Omega$
$$
\Prob(G) = p^{|E(G)|} (1-p)^{{n \choose 2} - |E(G)|} \,.
$$
This space may be viewed as ${n \choose 2}$ independent coin flips, one for each pair of vertices, where the probability of success (that is, drawing an edge) is equal to $p.$ Note that $p=p(n)$ can tend to zero as $n$ tends to infinity. All asymptotics throughout are as $n \rightarrow \infty $ (we emphasize that the notations $o(\cdot)$ and $O(\cdot)$ refer to functions of $n$ whose growth is bounded but are not necessarily positive). We say that an event in a probability space holds \emph{asymptotically almost surely} (\emph{a.a.s.}) if the probability that it holds tends to $1$ as $n$ goes to infinity. 

For $p\in (0,1)$ or $p=p(n) > 0$ tending to $0$ with $n$, define $\mathbb{L}n=\log_{\frac{1}{1-p}}n.$ For constant $p$, clearly $\mathbb{L}n=\Theta(\log n)$, but for $p=o(1)$ we have
$$
\mathbb{L}n= \frac {\log n}{-\log (1-p)} = (1+o(1)) \frac {\log n}{p}.
$$

\bigskip

Preliminary results for random graphs have been proved in~\cite{BCPWZ} where it is shown that for constant $p \in (0,1)$ and for $r < c \log_{\frac{1}{\min\{p,1-p\}}} n$ with $c < 1$, and also for constant $r$ and $p \gg n^{-1/r}$, $\sigma(G,r,m)=r-m+1$ a.a.s.\  (the required condition is that $pn^r \rightarrow \infty$). In this paper, we improve on these results using different techniques.
Our main results are summarized below. The proofs of the results in this subsection may be found in Section~\ref{sec:proofs}. Let us start with the following useful upper bound for the spy number.
\begin{theorem}\label{thm:upper_bound}
Let $\eps > 0$, $r=r(n) \in \N$, and $m=m(n) \in \N$. Consider a random graph $G \in G(n,p)$ with $p = p(n)< 1-\eps$. Then, a.a.s. 
$$
\sigma(G,r,m) \le \frac {r}{m} + 2(2+\sqrt{2}+\eps) \LL n.
$$
\end{theorem}
In particular, it follows from this theorem that $\sigma(G,r,m) = (1+o(1)) r/m$ whenever $r/m \gg \LL n$, and $\sigma(G,r,m) = \Theta(r/m)$ if $r/m = \Theta(\LL n)$. The next theorem provides a lower bound.
\begin{theorem}\label{thm:lower_bound}
Let $\eps > 0$, $\eta \in (0,1/3]$, $r=r(n) \in \N$, and $m=m(n) \in \N$. Consider a random graph $G \in G(n,p)$ with $p = n^{-1/3+\eta+o(1)}<1-\eps$. Then, a.a.s. 
$$
\sigma(G,r,m) \ge \min \{r-m, 2.99 \eta \LL n\}+1.
$$
\end{theorem}

Theorems~\ref{thm:upper_bound} and~\ref{thm:lower_bound} give an asymptotic behaviour (up to a constant factor) of the spy number of random graphs with average degree $n^{2/3+\eta+o(1)}$ for some $\eta \in (0,1/3]$. Therefore, we obtain the whole picture for such dense random graphs.

\begin{corollary}\label{cor:23_to_1}
Let $\eps > 0$, $\eta \in (0,1/3]$, $r=r(n) \in \N$, and $m=m(n) \in \N$. Consider a random graph $G \in G(n,p)$ with $p = n^{-1/3+\eta+o(1)}<1-\eps$. Then, a.a.s. 
$$
\sigma(G,r,m) = 
\begin{cases}
r-m+1 & \text{ if } r-m \le 2.99 \eta \LL n \\
\Theta(\LL n) & \text{ if } r-m > 2.99 \eta \LL n \text{ and } r/m = O(\LL n) \\
(1+o(1)) r/m & \text{ if } r/m \gg \LL n.
\end{cases}
$$
\end{corollary}

With a bit more effort and one additional idea, one can extend this (tight) result to slightly sparser graphs.
\begin{theorem}\label{thm:12_to_23}
Let $\eta \in (0,1/6]$, $r=r(n) \in \N$, and $m=m(n) \in \N$. Consider a random graph $G \in G(n,p)$ with $p = n^{-1/2+\eta+o(1)}$. Then, a.a.s. 
$$
\sigma(G,r,m) = 
\begin{cases}
r-m+1 & \text{ if } r-m =O(1) \\
(1+o(1)) (r-m) & \text{ if } r-m \gg 1 \text{ and } r-m \le 1.99 \eta \LL n \\
\Theta(\LL n) & \text{ if } r-m > 1.99 \eta \LL n \text{ and } r/m = O(\LL n) \\
(1+o(1)) r/m & \text{ if } r/m \gg \LL n.
\end{cases}
$$
\end{theorem}

For very sparse graphs (that is, graphs with average degree $n^{1/2-\eta+o(1)}$ for some $\eta \in [0,1/2]$) we are less successful. However, a partial progress has been made. First of all, we managed to investigate the case $r-m=O(1)$ for graphs with average degree at least $\log^3 n$ \emph{unless} the average degree is of order $\sqrt{n \log n}$. The sub-case $r=O(1)$ (and even $r$ growing with $n$ slowly) seems to be easy to deal with (the revolutionaries could go to a fixed set of vertices, say $[r]$; since a.a.s.\ all pairs of them have as many common neighbours as expected, and for each pair there is a common neighbour not adjacent to any spy, they can join themselves and create an unguarded meeting in the same way as in the other proofs of dense graphs). Unfortunately, no approach to solve the general case $r-m=O(1)$ and $r$ tending (fast) to infinity with $n$ (and so $m$ as well) is known. It is not clear if this peculiar gap at $\sqrt{n \log n}$ is an outcome of a wrong approach used or perhaps the behaviour of the spy number changes in this window. 
This remains to be investigated.

\begin{theorem}\label{thm:sparse_graphs_constant}
Let $\eps > 0$, $r=r(n) \in \N$, and $m=m(n) \in \N$ such that $r-m = O(1)$. Consider a random graph $G \in G(n,p)$ with $\log^3 n \le pn \ll \sqrt{n \log n}$ or $\sqrt{n \log n} \ll pn < (1-\eps)n$. Then, a.a.s. 
$$
\sigma(G,r,m) = r-m+1.
$$
\end{theorem}

For $r-m \gg 1$ and the average degree of $n^{\eta+o(1)}$ for some $\eta \in (0,1/2)$ we have the following result. 

\begin{theorem}\label{thm:sparse_graphs_nonconstant}
Let $np=n^{\eta+o(1)}$ for some $\eta \in (0, \frac12)$. Let $\omega$ be any function tending to infinity arbitrarily slowly and suppose $r-m \gg 1$. Then, a.a.s. $$
\sigma(G,r,m) = (1+o(1))(r-m),
$$
provided that one of the following situations occurs:
\begin{enumerate}
\item [(i)] $r-m =o(\min \{pn, \frac{n}{(pn)^2}\})$ and $pn \geq (\omega n \log n)^{1/4}$,
\item [(ii)]  $r-m =o(pn/(\omega \log n))$ and $(\omega n \log n)^{1/6} \leq pn < (\omega n \log n)^{1/4}$,
\item [(iii)] $r-m =o(pn)$ and $pn < (\omega n \log n)^{1/6}$.
\end{enumerate}
\end{theorem}
Let $f(\eta) = 1-\eta$ and let
$$ 
g(\eta) = 
\begin{cases}
\eta & \text{ if } \eta \in (0, 1/3] \\
1-2\eta & \text{ if } \eta \in (1/3,1/2).
\end{cases}
$$
Note that for a given $\eta \in (0,1/2)$, Theorems~\ref{thm:sparse_graphs_constant} and~\ref{thm:sparse_graphs_nonconstant} investigate the case $r-m = n^{\alpha+o(1)}$ for $\alpha \in [0, g(\eta))$; the case $r-m = n^{\alpha+o(1)}$ for $\alpha > f(\eta)$ follows from Theorem~\ref{thm:upper_bound}. The behaviour of the spy number when $r-m = n^{\alpha+o(1)}$ for some $\alpha \in [g(\eta),f(\eta)]$ is not known and remains to be analyzed.

\bigskip

Over the last few years there was an explosion of research related to  the game of \emph{Cops and Robbers}, introduced independently by Nowa\-kowski and Winkler~\cite{nw} and Quilliot~\cite{q} almost thirty years ago. While much of the earlier work focused on deterministic graphs, a number of papers on random graphs generated a lot of interest recently~\cite{bkl,bpw,lp2,p,pw,pw2}. Is the game of Revolutionaries and Spies going to follow the same path? For a survey of results on vertex pursuit games such as Cops and Robbers, the reader is directed to the surveys~\cite{al,fomin,gena} and the recent monograph~\cite{AB1}.

\section{Proofs}\label{sec:proofs}

Let $\Ss(v,i)$ denote the set of vertices whose distance from $v$ is precisely $i$, and $\Nn(v,i)$ the set of vertices whose distance from $v$ is at most $i$. (In particular, $\Nn(v)=\Ss(v,1)$, the neighbourhood of $v$.) Also, $\Nn[S,i] = \bigcup_{v \in S} \Nn(v,i)$, $\Nn[S]$ denotes $\Nn[S,1]$, the closed neighbourhood of $S$, and $\Nn(S) = \Nn[S] \setminus S$, the (open) neighbourhood of $S$.  Finally, $\Nn^c(v)$ denotes $V(G) \setminus \Nn(v,1)$, the non-neighbourhood of $v$.

\subsection{Proof of Theorem~\ref{thm:upper_bound}}

Let us start by showing some typical properties of a random graph $G(n,p)$. First, we calculate a universal upper bound for the size of an intersection of non-neighbourhoods. 
\begin{lemma}\label{lem:non-neigh}
Let $\alpha > 0$, $\beta > \frac {1+\alpha}{\alpha}$, and $\eps > 0$. Consider a random graph $G \in G(n,p)$ with $p < 1-\eps$. Then, a.a.s.\ for every set $S \subseteq V(G)$ of cardinality $\beta \LL n$ we have
$$
\left| \bigcap_{v \in S} \Nn^c(v) \right| \le \alpha \beta \LL n.
$$
\end{lemma}

From this lemma we get immediately the following corollary.

\begin{corollary}\label{cor:non-neigh}
Let $\beta > 1$ and $\eps > 0$. Consider a random graph $G \in G(n,p)$ with $p < 1-\eps$. Then, a.a.s.\ for every set $S \subseteq V(G)$ of cardinality $\beta \LL n$ we have
$$
\left| \bigcap_{v \in S} \Nn^c(v) \right| \le (1+\eps) \frac {\beta}{\beta-1} \LL n.
$$
\end{corollary}

\begin{proof}[Proof of Lemma~\ref{lem:non-neigh}]
Let $\alpha > 0$ and $\beta > \frac {1+\alpha}{\alpha}$. If $\alpha \beta \LL n \ge n$, then there is nothing to prove; the bound trivially holds. Suppose then that $\alpha \beta \LL n < n$. Let $X$ be the number of pairs of sets $(S,T)$ such that $S,T \subseteq V(G)$, $s=|S|=\beta \LL n$, $|T|=\alpha |S|$, and there is no edge between $S$ and $T$. Since $(1-p)^s = n^{-\beta}$, we get that
\begin{eqnarray*}
\E [X] &\le& {n \choose s} {n \choose \alpha s} (1-p)^{\alpha s^2} \\
&\le& n^{s + \alpha s} n^{- \alpha \beta s} \\
&=& \exp \left( (1+\alpha) \beta (\LL n) (\log n) - \alpha \beta^2 (\LL n) (\log n) \right) \\
&=& \exp \left( (1+\alpha-\alpha \beta) \beta (\LL n) (\log n) \right) \\
&=& o(1).
\end{eqnarray*}
(Note that $1 + \alpha - \alpha \beta < 0$.) The result follows from Markov's inequality.
\end{proof}

It is not difficult to show that a.a.s.\ a random graph $G(n,p)$ has a dominating set of size $(1+o(1)) \LL n$. Now, we will show that a slightly larger set can not only dominate the rest of the graph but also can create a matching with any set of cardinality $O( \LL n)$.

\begin{lemma}\label{lem:matching}
Let $\eps>0$ and consider a random graph $G \in G(n,p)$ with $p = n^{-\eta+o(1)}<1-\eps$ for some $\eta \in [0,1]$. Let $\gamma, \delta > 0$ be such that $\gamma-\delta > 1 + \eta$.  Then, a.a.s.\ there exists a set $S \subseteq V(G)$ of cardinality $\gamma \LL n$ such that for all sets $T \subseteq V(G) \setminus S$ of size at most $\delta \LL n$ there is a perfect matching from $T$ to some subset $S'$ of $S$.
\end{lemma}

\begin{proof}
Let $\gamma, \delta > 0$ be such that $\gamma-\delta > 1 + \eta$. If $\gamma \LL n \ge n$, then there is nothing to prove; the claim trivially holds. Suppose then that $\gamma \LL n < n$. We will use Hall's theorem for perfect matchings in bipartite graphs to prove the result. We need to show that a.a.s.\ Hall's condition is satisfied for any set $T \subseteq V(G) \setminus S$ of size at most $\delta \LL n$, that is, we have to show that $|\Nn(T) \cap S| \ge |T|$. For $1 \le t \le \delta \LL n$, let $X_t$ be the random variable counting the number of sets $T$ of cardinality $t$ for which the condition fails. For a given $t$ in this range, we get
\begin{eqnarray*}
\E [X_t] &\le& {n \choose t} {\gamma \LL n \choose t-1} (1-p)^{t (\gamma \LL n - (t-1))},
\end{eqnarray*}
since there are at most ${n \choose t}$ ways to choose $T$, ${\gamma \LL n \choose t-1}$ ways to choose vertices in $S$ that are possibly connected to $T$, and all other vertices of $S$ are not adjacent to any vertex of $T$. We have
\begin{eqnarray*}
\E [X_t] &\le& \exp \big( t \log n + t \log (\gamma \LL n) \big) (1-p)^{t (\gamma - \delta) \LL n} \\
&=& \exp \Big(  (1+o(1)) \big( 1 + \eta \big) t \log n \Big) n^{-t (\gamma - \delta)} \\
&=& \exp \Big(  (1+o(1)) \big( 1 + \eta -\gamma + \delta \big) t \log n \Big) \\
&=& \exp \Big(  - (1+o(1)) 2 \eps t \log n \Big) \\
&\le& \exp \Big(  - \eps t \log n \Big),
\end{eqnarray*}
where $\eps = (\gamma-\delta-1-\eta)/2 > 0$. Finally, 
$$
\E \left[ \sum_{t=1}^{\delta \LL n} X_t  \right] = O(\exp(-\eps \log n)) = o(1)
$$ 
so a.a.s.\ Hall's condition fails  for no set $T$ under consideration by Markov's inequality. The proof is finished.
\end{proof}

Now, we are ready to prove the main result of this subsection.

\begin{proof}[Proof of Theorem~\ref{thm:upper_bound}]
Consider a random graph $G \in G(n,p)$ with $p = n^{-\eta+o(1)} < 1-\eps$ for some $\eta \in [0,1]$. Let
$$
\delta = \frac {1-\eta+\sqrt{\eta^2 + 2\eta + 5}}{2} > 1
$$
and 
$$
\gamma = 1+\eta+\delta+\eps = \frac {3+\eta+\sqrt{\eta^2 + 2\eta + 5}}{2} + \eps \le 2 + \sqrt{2} + \eps.
$$
It follows from Lemma~\ref{lem:matching} that a.a.s.\ there exists a set $A \subseteq V(G)$ of cardinality $\gamma \LL n$ such that for all sets $T \subseteq V(G) \setminus A$ of size at most $\delta \LL n$ there is a perfect matching from $T$ to some set $A' \subseteq A$. 

We will split spies into three groups: the first two groups, $\tau_1$ and  $\tau_2$, consist each of $\gamma \LL n$ \emph{super-spies}, and the third group $\tau_3$ consists of $\lfloor r/m \rfloor$ \emph{regular-spies}. Super-spies from team $\tau_1$ will occupy the whole set $A$ at odd times but some of them might be sent to a mission at even times. If this is the case, then they will be back to the set $A$ in the next round. Similarly, super-spies from team $\tau_2$ will occupy $A$ at even times but might be used to protect some other vertices at odd times. In particular, the set $A$ will be constantly protected and so no unguarded meeting can take place there. Regular-spies (team $\tau_3$) will occupy a set $B_t \subseteq V(G) \setminus A$ at time $t$. Moreover, no two regular-spies will occupy the same vertex, so $|B_t| = \lfloor r/m \rfloor$ for all $t$.

The revolutionaries start the game by occupying a set $R_1$ and forming at most $\lfloor r/m \rfloor$ meetings. The super-spies (from both teams $\tau_1$ and $\tau_2$) go to the set $A$. The regular-spies can easily protect the vertices of $V(G) \setminus A$ in which meetings take place by placing a spy on each vertex where there are at least $m$ revolutionaries. The remaining regular-spies go to arbitrary vertices of $V(G) \setminus A$ not occupied by another spy. As a result, no unguarded meeting takes place in the first round.

Suppose that no unguarded meeting takes place at time $t-1$ and that the regular-spies occupy a set $B_{t-1} \subseteq V(G) \setminus A$ of cardinality $\lfloor r/m \rfloor$. At the beginning of round $t$, the revolutionaries might form at most $\lfloor r/m \rfloor$ meetings at vertices of $M_t \subseteq V(G) \setminus A$ (as we already pointed out, meetings at $A$ are constantly protected by super-spies, so we do not have to worry about them). Let $B = M_t \cap B_{t-1}$ be the set of vertices in which meetings are already guarded by regular-spies. It remains to show that there exists a perfect matching between $M_t \setminus B$ and some subset $S$ of $A \cup (B_{t-1} \setminus B) $. Indeed, if this is the case, then the regular-spies that do not protect any meeting as well as the super-spies from the team protecting $A$ in the previous round, move from $S$ to $M_t \setminus B$. The Super-spies from another team come back to $A$ to guard this set and prepare themselves for another mission. No unguarded meeting takes place at time $t$, $B_t \subseteq V(G) \setminus A$, and no two regular-spies occupy the same vertex. The result will follow by induction on $t$.

In order to show that a perfect matching can be formed we use Hall's theorem. For a given $T \subseteq M_t \setminus B$,  we need to show that $|\Nn(T) \cap (A \cup (B_{t-1} \setminus B))| \ge |T|$. If $|T| \le \delta \LL n$, then a perfect matching from $T$ to some $A' \subseteq A$ exists and so 
$$
|\Nn(T) \cap (A \cup (B_{t-1} \setminus B))| \ge |\Nn(T) \cap A| \ge |T|.
$$
Hall's condition then holds for such small sets. Suppose then that $|T| > \delta \LL n$. It follows from Corollary~\ref{cor:non-neigh} that all but at most 
\begin{eqnarray*}
\left( 1+ \frac {\eps}{4} \right) \frac {\delta}{\delta -1} \LL n &=& \left( 1+ \frac {\eps}{4} \right) \frac {3+\eta+\sqrt{\eta^2 + 2\eta + 5}}{2} \LL n \\
&<& \left( \frac {3+\eta+\sqrt{\eta^2 + 2\eta + 5}}{2} + \eps \right) \LL n\\
&=& \gamma \LL n
\end{eqnarray*}
vertices of $V(G) \setminus T$ have at least one neighbour in $T$. Hence,
\begin{eqnarray*}
|\Nn(T) \cap (A \cup (B_{t-1} \setminus B))| &>& |A| + \lfloor r/m \rfloor - |B| - \gamma \LL n \\
&=& \lfloor r/m \rfloor - |B| \\
&\ge& |M_t| - |B| \ge |T|.
\end{eqnarray*}
Thus, Hall's condition holds for larger sets as well and the proof is finished.
\end{proof}

\subsection{Proof of Theorem~\ref{thm:lower_bound}}

For the proofs of the lower bounds in Theorem~\ref{thm:lower_bound}, we employ the following adjacency property and its generalizations. For fixed positive integers $k$ and $l$, we say that a graph $G$ is $(l,k)$-\emph{existentially closed} (or $(l,k)$-\emph{e.c.}) if for any two disjoint subsets of $V(G)$, $A \subseteq V(G), B \subseteq V(G)$, with $|A|=l$ and $|B|=k$, there exists a vertex $z \in V(G) \setminus (A \cup B)$ not joined to any vertex in $B$ and joined to every vertex in $A$. We will use the following simple observation.

\begin{theorem}\label{thm:k_ec}
Let $r, m, s$ be positive integers such that $s \le r-m$, and let $G$ be any $(2,s)$-e.c.\ graph. Then
$$
\sigma(G,r,m) \ge s+1.
$$
In particular, if $G$ is $(2,r-m)$-e.c., then $\sigma(G,r,m)=r-m+1$.
\end{theorem}
\begin{proof}
Suppose that $r$ revolutionaries play the game on a graph $G$ against $s$ spies. The revolutionaries start by occupying $r$ distinct vertices. No matter what $s$ spies do, there will be $r-s \ge m$ unguarded revolutionaries. Since $G$ is $(2,s)$-e.c., any two of them can meet in the next round and stay unguarded. In the following round, another revolutionary can join the two forming a group of three unguarded revolutionaries. This argument can be repeated (each time at least one more revolutionary joins the group) until an unguarded meeting of $m$ revolutionaries is formed, and the game ends. The result holds.
\end{proof}

It remains to investigate for which values of $s$ a random graph is $(2,s)$-e.c.\ a.a.s. Since we would like to match an upper bound of $O(\LL n)$, our goal is to obtain a lower bound of $\Omega(\LL n)$. Hence, the graph must be dense enough.

\begin{lemma}\label{lem:k_ec-gnp}
Consider a random graph $G \in G(n,p)$ with $p = n^{-1/3+\eta+o(1)}<1-\eps$ for some $\eta \in (0,1/3]$ and $\eps>0$. Then, for
$$
s = 2.99 \eta \LL n
$$ 
we have that a.a.s.\ $G$ is $(2,s)$-e.c.
\end{lemma}
\begin{proof}
Fix any two disjoint subsets of $V(G)$, $A, B$, with $|A|=2$ and $|B|=s$. For a vertex $z \in V(G) \setminus (A \cup B)$, the probability that $z$ is joined to both vertices of $A$ and no vertex of $B$ is $p^2(1-p)^s$. Since the edges are chosen independently, the probability that no suitable vertex can be found for this particular choice of $A$ and $B$ is $(1-p^2(1-p)^s)^{n-s-2}$. 

Let $X$ be the random variable counting the number of pairs of $A,B$ for which no suitable $z$ can be found. We have 
\begin{eqnarray*}
\E [X] &=& {n \choose 2}{n-2 \choose s} (1-p^2(1-p)^s)^{n-s-2} \\
&\le& n^{s+2} \exp \left( - p^2(1-p)^s (n-s-2) \right) \\
&=& \exp \left( (s+2) \log n - p^2(1-p)^s n (1+o(1)) \right).
\end{eqnarray*}
If $p=\Theta(1)$, then $s = \Theta(\log n)$ and so
\begin{eqnarray*}
\E [X] &\le& \exp \left( O(\log^2 n) - \Omega(n^{1-2.99\eta}) \right) \\
&=& \exp \left( O(\log^2 n) - \Omega(n^{0.003}) \right) \\
&=& o(1).
\end{eqnarray*}
For $p = o(1)$ we have $s = (1+o(1)) 2.99 \eta (\log n) / p = n^{1/3-\eta+o(1)}$ and so
\begin{eqnarray*}
\E [X] &\le& \exp \left( n^{1/3-\eta+o(1)} - n^{2(-1/3+\eta+o(1)) - 2.99 \eta + 1} \right)\\
&=& \exp \left( n^{1/3-\eta+o(1)} - n^{1/3-0.99 \eta+o(1)} \right) \\
&=& o(1).
\end{eqnarray*}
The result follows by Markov's inequality.
\end{proof}

Theorem~\ref{thm:lower_bound} follows immediately from Theorem~\ref{thm:k_ec} and Lemma~\ref{lem:k_ec-gnp}.

\begin{proof}[Proof of Theorem~\ref{thm:lower_bound}]
It follows from Lemma~\ref{lem:k_ec-gnp} that a.a.s.\ $G$ is $(2,s_0)$-e.c.\ for $s_0=2.99 \eta \LL n$. If $r-m \ge s_0$, then we get that $\sigma(G,r,m) \ge s_0+1$ by Theorem~\ref{thm:k_ec} (applied with $s=s_0$). Suppose then that $r-m<s_0$. Note that Theorem~\ref{thm:k_ec} cannot be applied with $s=s_0$ anymore. However, it is clear that any $(2,s_0)$-e.c.\ graph is also (deterministically) $(2,r-m)$-e.c. Using Theorem~\ref{thm:k_ec} again (this time with $s=r-m$) we get that $\sigma(G,r,m) = r-m+1$, and the proof is finished.
\end{proof}

\subsection{Proof of Theorem~\ref{thm:12_to_23}}

In this subsection, we investigate slightly sparser graphs than before, namely, graphs with average degree tending to infinity faster than $\sqrt{n \log n}$. Let $G \in G(n,p)$ be a random graph with $\sqrt{\log n / n} \ll p < n^{-0.33}$. (The previous results cover graphs with average degree at least $n^{2/3+\eta+o(1)}$ for some $\eta>0$. Hence, we may assume that the average degree is at most $n^{0.67}$.) Let us first note that $G$ is a.a.s.\ $(2,s)$-e.c., provided that $s=O(1)$. Indeed, the probability that this is not the case can be bounded from above by
$$
n^{s+2} \left( 1 - p^2 (1-p)^s \right)^{n-o(n)} \le \exp  \left( O(s \log n) - \Omega(p^2 (1-p)^s n) \right),
$$
which tends to zero even for $s=s(n)$ going to infinity with $n$ slow enough. Therefore, it follows from Theorem~\ref{thm:k_ec} that a.a.s.\ $\sigma(G,r,m)=r-m+1$, provided that $r-m=O(1)$.

In order to deal with $r-m$ tending to infinity with $n$ faster than before, we need to relax slightly the $(2,s)$-e.c.\ property at the price of obtaining a bit weaker lower bound for the spy number. For fixed positive integers $l$ and $k$, we say that a graph $G$ is $(1,l,k)$-\emph{existentially closed} (or $(1,l,k)$-\emph{e.c.}) if for any vertex $v \in V(G)$ and two disjoint subsets of $A, B \subseteq V(G) \setminus \{v \}$ with $|A|=l$ and $|B|=k$, there exists a vertex $z \in V(G) \setminus (\{v \} \cup A \cup B)$ not joined to any vertex in $B$, joined to $v$, and joined to \emph{some} vertex in $A$. Note that $(1,l,k)$-e.c.\ is a natural generalization of the $(2,k)$-e.c.\ property, which is equivalent to $(1,1,k)$-e.c. Any $(1,l_1,k)$-e.c.\ graph is also $(1,l_2,k)$-e.c., provided that $l_1 < l_2$. Moreover, any $(1,l,k_1)$-e.c.\ graph is also $(1,l,k_2)$-e.c., provided that $k_1 > k_2$. We will use the following observation (since this is a simple generalization of Theorem~\ref{thm:k_ec}, the proof of it is omitted).

\begin{theorem}\label{thm:lk_ec}
Let $r, m, s, l$ be positive integers such that $s \le r-m-l+1$, and let $G$ be any $(1,l,s)$-e.c.\ graph. Then
$$
\sigma(G,r,m) \ge s+1.
$$
\end{theorem}

As before, let $G \in G(n,p)$ with $\sqrt{\log n / n} \ll p < n^{-0.33}$. Suppose first that $ps=o(1)$. Take $\eps = \eps(n) = 3 \log n / (p^2 n) = o(1)$ and $l=\eps s$. The probability that $G$ is not $(1,l,s)$-e.c.\ is bounded from above by
\begin{eqnarray*}
n^{s+\eps s +1} \Big( 1 - p (1-(1-p)^{\eps s}) (1-p)^s \Big)^{n-o(n)} &\le& \exp  \Big( 2s \log n - (1+o(1)) p (p\eps s) n \Big) \\
&=& \exp  \Big( s ( 2 \log n - (1+o(1)) p^2 \eps n) \Big) \\
&=& \exp  \Big( - (1+o(1)) s \log n  \Big) \\
&=& o(1).
\end{eqnarray*}
Hence, $G$ is $(1,\eps s, s)$-e.c.\ a.a.s., provided that $ps=o(1)$. After applying Theorem~\ref{thm:lk_ec} to $s=(r-m)/(1+\eps)$ and $l=\eps s = \eps (r-m)/(1+\eps)$, we obtain that a.a.s.\ $\sigma(G,r,m) \ge s = (1+o(1)) (r-m)$, provided that $r-m=o(1/p)$. Since $r-m+1$ is a (deterministic) trivial upper bound for any graph, the conclusion is that a.a.s.\ $\sigma(G,r,m) = (1+o(1)) (r-m)$.

Suppose now that $ps=c \in \R$. This time we get a bound of 
\begin{eqnarray*}
n^{s+\eps s +1} \Big( 1 - p (1-(1-p)^{\eps s}) (1-p)^s \Big)^{n-o(n)} &\le& \exp  \Big( 2s \log n - (1+o(1)) p (1-e^{-c\eps}) e^{-c} n \Big) \\
&=& \exp \Big( 2 c \log n / p - (1+o(1)) p (1-e^{-c\eps}) e^{-c} n \Big) \\
&=& \exp \Big( O(\log n / p) - \Omega( p (1-e^{-c \eps}) n) \Big),
\end{eqnarray*}
which tends to zero for $\eps$ tending to zero slow enough (recall that $p \gg \sqrt{\log n / n}$). As a result, $G$ is $(1,\eps s, s)$-e.c.\ a.a.s., provided that $s = \Theta(1/p)$. As before, it follows from Theorem~\ref{thm:lk_ec} that a.a.s.\ $\sigma(G,r,m) = (1+o(1)) (r-m)$, provided that $r-m=\Theta(1/p)$.

Finally, suppose that $ps \gg 1$ but $s \le c \LL n$ (the constant $c$ will be determined soon). This time we assume that the average degree satisfies $np = n^{1/2+\eta+o(1)}$ for some $\eta \in (0,1/6]$ and that $\eps = 1/(ps) = o(1)$ to get a bound of 
\begin{eqnarray*}
n^{s+\eps s +1} \Big( 1 - p (1-(1-p)^{\eps s}) (1-p)^s \Big)^{n-o(n)} &\le& \exp  \Big( 2s \log n - (1+o(1)) p (1-e^{-1}) (1-p)^s n \Big) \\
&=& \exp \Big( O(\log^2 n / p) - \Omega( p n^{1-c} ) \Big) \\
&=& \exp \Big( n^{1/2-\eta+o(1)} - \Omega( n^{1/2-c+\eta} ) \Big)
\end{eqnarray*}
which tends to zero for, say, $c=1.99 \eta$. We have that $G$ is $(1,\eps s, s)$-e.c.\ a.a.s., and so a.a.s.\ $\sigma(G,r,m) = (1+o(1)) (r-m)$, provided that $r-m \le 1.99 \eta \LL n$. For $r-m > 1.99 \eta \LL n$ we apply Theorem~\ref{thm:lk_ec} with $s = 1.99 \eta \LL n / (1+\eps)$ and $l = \eps s$ to get an asymptotically almost sure lower bound of $(1+o(1)) 1.99 \eta \LL n$. This finishes the proof of Theorem~\ref{thm:12_to_23}, since the upper bounds in the last two cases follow directly from Theorem~\ref{thm:upper_bound}.

\subsection{Proof of Theorem~\ref{thm:sparse_graphs_constant}}

Let us start with the following expansion-type properties of random graphs. 

\begin{lemma}\label{lem:gnp exp}
Suppose that $d=p(n-1) \ge \log^3 n$. Let $G \in G(n,p)$.  The following property holds a.a.s. Let $S \subseteq V(G)$ be any set of $s=|S|$ vertices, and let $i \in \N$. If $s$ and $i$ are such that $s d^i < n / \log n$, then
$$
\left| \Nn [S,i] \right| = (1+o(1)) s d^i.
$$
\end{lemma}

\begin{proof} 
Let $S \subseteq V(G)$, $s=|S|$, and consider the random variable $X = X(S) = |\Nn[S]|$. We will bound $X$ in a stochastic sense. There are two things that need to be estimated: the expected value of $X$, and the concentration of $X$ around its expectation. 

It is clear that
\begin{eqnarray*}
\E [X] &=& n - \left(1- \frac {d}{n-1} \right)^s (n-s) \\
&=& n - \exp \left( - \frac {ds}{n} (1+O(d/n)) \right) (n-s) \\
&=& ds (1+O(\log^{-1} n)) 
\end{eqnarray*}
provided $ds \le n/ \log n $. We next use a consequence of Chernoff's bound (see e.g.~\cite[p.\ 27, Corollary~2.3]{JLR}), that 
\begin{equation}\label{chern}
\Prob( |X-\E X| \ge \eps \E X) ) \le 2\exp \left( - \frac {\eps^2 \E X}{3} \right)  
\end{equation}
for  $0 < \eps < 3/2$. This implies that the expected number of sets $S$ that have  $\big| |\Nn[S]| - d|S| \big| > \eps d|S|$ and $|S| \le n/(d\log n)$  is, for $\eps = 2/{\log n}$, at most
$$
\sum_{s \ge 1} 2 n^s \exp \left( - \frac {\eps^2 s \log^3 n}{3+o(1)} \right) =  \sum_{s \ge 1} 2 \exp \left( s \log n - \left( \frac 43 +o(1) \right) s \log n \right)=o(1).
$$
So a.a.s.\ if $ |S| \le n/d\log n\mbox{ then } |\Nn[S]| =d |S| (1+O(1/\log n))$ where the bound in $O()$ is uniform. We may assume  this statement holds deterministically. 
 
Given this assumption, we have good bounds on the ratios of the cardinalities of $\Nn[S]$, $\Nn[\Nn[S]] = \Nn[S,2]$,  and so on. We consider this up to the $i$'th iterated neighbourhood  provided $sd^i \le  {n}/{\log n}$ and thus  $i = O(\log n /\log \log n)$. Then the cumulative multiplicative error term  is $(1+O(\log^{-1}n))^i = (1+o(1))$, that is, 
$$
|\Nn[S,i]| = (1+o(1)) sd^{i} 
$$
for all $s$ and $i$ such that $sd^i \le n / \log n$. The proof is finished.
\end{proof}

Here is another technical lemma that will be used.

\begin{lemma}\label{lem:gnp exp2}
Suppose that $d=p(n-1) \ge \log^3 n$. Let $G\in G(n,p)$. The following properties hold a.a.s.

\begin{enumerate}
\item [(i)] Let $S \subseteq V(G)$ be any set of $s=|S|$ vertices, $x \in V(G) \setminus S$, $y \in V(G) \setminus (S \cup \{x\})$, and let $i \in \N$. If $s$ and $i$ are such that $s = O(1)$ and $(s+2) d^i < n / \log n$, then
\begin{eqnarray}
\left| \Nn(x,i) \setminus \Nn[S,i] \right| &=& (1+o(1)) d^i \label{eq:cond_x} \\
\left| \Nn(y,i) \setminus \Nn[S \cup \{x\},i] \right| &=& (1+o(1)) d^i \label{eq:cond_y}
\end{eqnarray}
\item [(ii)] Let $S \subseteq V(G)$ be any set of $s=|S|$ vertices, $x \in V(G) \setminus S$, and let $i \in \N$. Suppose that $d = n^{\eta+o(1)}$ for some $\eta \in (0,1/2)$. If $s$ and $i$ are such that  $i \eta < 1-\eta$ and $1.01 i s \le (1-\eta -i \eta) d/4$, then
$$
\left| \Nn(x,i) \setminus \Nn[S,i] \right| \ge d^i / 2.
$$
\end{enumerate}
\end{lemma}
\begin{proof}
The proof of part (i) is easy. Since we aim for a statement that holds a.a.s., we can assume that the property stated in Lemma~\ref{lem:gnp exp} holds deterministically.  We have $\left| \Nn[S,i] \right|=(1+o(1)) sd^i$ and $\left| \Nn[S \cup \{x\},i] \right|=(1+o(1)) (s+1)d^i$, which immediately implies that~(\ref{eq:cond_x}) holds. Since $\left| \Nn[S \cup \{x,y\},i] \right|=(1+o(1)) (s+2)d^i$, exactly the same argument can be used to show~(\ref{eq:cond_y}). This establishes part (i).

\smallskip

For part (ii), for simplicity set $\eps=1/2 - \eta$. We use a well known fact that a random graph a.a.s.\ does not contain a copy of a dense graph as a subgraph. Formally, a.a.s.\ there is no finite subgraph $\hat{G}$ on $\hat{n}$ vertices and $\hat{m}$ edges, provided that
$$
n^{\hat{n}} p^{\hat{m}} = n^{\hat{n}} n^{(-1+\eta+o(1))\hat{m}} = n^{\hat{n} + \hat{m} (-1/2-\eps+o(1))} \to 0.
$$
(This is, in fact, an easy observation, since the expected number of copies of $\hat{G}$ is of order $n^{\hat{n}} p^{\hat{m}}$ and the claim holds by Markov's inequality.) Therefore, a sufficient condition for not having $\hat{G}$ as a subgraph is that $\hat{n}/\hat{m} < 1/2+\eps$. In particular, a.a.s.\ there is no $K_{2,\lfloor 1/\eps \rfloor+1}$ as a subgraph, since 
$$
(2+\lfloor 1/\eps \rfloor+1)/(2 (\lfloor 1/\eps \rfloor+1)) = 1/2 + 1/(\lfloor 1/\eps \rfloor+1) < 1/2 + \eps. 
$$
Moreover, a.a.s.\ no two adjacent vertices share $\lfloor 1/\eps \rfloor$ or more neighbours. (Note that this structure has one vertex less and one edge less than $K_{2,\lfloor 1/\eps \rfloor+1}$, and so it is denser.) Therefore, we may assume that for any two distinct vertices $u$ and $w$, 
\begin{equation}\label{eq:lose}
|\Ss(u,1) \cap \Nn(w,1)| \le 1/\eps.
\end{equation}

Now, let us come back to the proof of part (ii). In order to illustrate the main idea, let us consider the case $i=1$ first. It follows from Lemma~\ref{lem:gnp exp} that $|\Nn(x,1)|=(1+o(1))d$ and that $\left| \Nn[S,1] \right|=(1+o(1)) sd$. On the other hand,~(\ref{eq:lose}) implies that the number of common neighbours of $S$ and $x$ is at most $s/\eps$. Hence, there are at least $(1+o(1))(d-s/\eps)$ neighbours of $x$ that are not adjacent to any vertex of $S$.  We will keep using (inductively) a similar idea to bound from below the number of vertices at distance $i$ from $x$ that are at distance at least $i+1$ from $S$. Let us call such a set of vertices $i$-\emph{good}. Unfortunately,~(\ref{eq:lose}) is no longer good enough to do it, since, for example, a given neighbour of some vertices of $S$ can potentially share a few neighbours with all neighbours of $x$. In order to rule this out we need to consider a tree rooted at vertex $x$ to obtain a denser forbidden subgraph. 

Suppose that a vertex $u$ at distance $i-1$ from $S$ shares $t$ neighbours with some vertices at distance $i-1$ from $x$. Consider a subgraph consisting of paths of lengths $i$ from $x$ to these neighbours together with $t$ edges connecting $u$ with them. In the extreme case, we get a graph on $\hat{n} = ti+2$ vertices and $\hat{m}=ti+t$ edges. (All other possibilities yield denser graphs.) We obtain a forbidden subgraph provided that
$$
t > \frac {2}{1/2+\eps - i (1/2 - \eps)} =: c = c(i,\eps).
$$
(Recall that $i (1/2 - \eps) < 1/2+\eps$.) Hence, we may assume that the number of shared neighbours is at most $\lfloor c \rfloor \le c$.

The rest of the proof is straightforward. Suppose that a lower bound of $d^{i-1} (1 - 1.01(i-1)cs/d)$ for the size of an $(i-1)$-good set is obtained. Lemma~\ref{lem:gnp exp} implies that this set expands by a factor of $(1+o(1))d$, but only a subset of it is $i$-good. Moreover, it follows from the same lemma that $\left| \Nn[S,i-1] \right|=(1+o(1)) sd^{i-1}$ and so at most $(1+o(1))csd^{i-1}$ vertices are eliminated and the remaining ones form an $i$-good set. We obtain a lower bound for the size of an $i$-good set of 
$$
d^{i} (1 - 1.01 (i-1)cs/d) - (1+o(1))csd^{i-1} \ge d^i (1 - 1.01ics/d).
$$
Since $1.01 i s \le (1/2+\eps-i(1/2-\eps)) d/4$, the bound is at least $d^i/2$ and the assertion holds.
\end{proof}

Now, we are almost ready to prove Theorem~\ref{thm:sparse_graphs_constant}. However, the $(l,k)$-e.c.\ property needs to be generalized again. For fixed positive integers $j$, $k$ and $l$, we say that a graph $G$ is $(l,k)_j$-\emph{existentially closed} (or $(l,k)_j$-\emph{e.c.}) if for any two disjoint subsets of $V(G)$, $A \subseteq V(G), B \subseteq V(G)$, with $|A|=l$ and $|B|=k$, there exists a vertex $z \in V(G) \setminus (A \cup B)$ at distance at most $j$ from every vertex in $A$ and at distance at least $j+1$ from every  vertex in $B$. This is indeed a generalization of the $(l,k)$-e.c.\ property, which is equivalent to $(l,k)_1$-e.c. Theorem~\ref{thm:k_ec} can be generalized easily to this new definition.

\begin{theorem}\label{thm:k_ec_i_distance}
Let $j, r, m, s$ be positive integers such that $s \le r-m$, and let $G$ be any $(2,s)_j$-e.c.\ graph. Then
$$
\sigma(G,r,m) \ge s+1.
$$
In particular, if $G$ is $(2,r-m)_j$-e.c., then $\sigma(G,r,m)=r-m+1$.
\end{theorem}

\begin{proof}[Proof of Theorem~\ref{thm:sparse_graphs_constant}]
The statement for dense random graphs (that is, graphs with average degree satisfying $\sqrt{n \log n} \ll pn < (1-\eps)n$) is already proved---see Corollary~\ref{cor:23_to_1}, Theorem~\ref{thm:12_to_23}, and the proof of Theorem~\ref{thm:12_to_23} for the case $pn = n^{1/2+o(1)}$ in the previous subsection. It remains to investigate sparse random graphs (that is, graphs with average degree satisfying $\log^3 n \le pn \ll \sqrt{n \log n}$).  The proof will follow from Theorem~\ref{thm:k_ec_i_distance} once it is shown that there exists an integer $j$ such that for all $s=O(1)$ a.a.s.\ the graph is $(2,s)_j$-e.c. 

Let us first consider $p \ge c/\sqrt{n}$ for some $c > 0$ and $p \ll \sqrt{\log n / n}$. For this range of $p$ we will use $j=2$. Let $S \subseteq V(G)$ ($s=|S|$), $x \in V(G) \setminus S$, and $y \in V(G) \setminus (S \cup \{x\})$. Lemma~\ref{lem:gnp exp} and Lemma~\ref{lem:gnp exp2}(i) with $i=1$ can be used to get that the size of neighbourhood of $S$ is $(1+o(1))spn$, and that there are two disjoint sets $X$ and $Y$ that are also disjoint from $\Nn[S]$, each of cardinality $(1+o(1))pn$, such that $X \subseteq \Nn(x)$ and $Y \subseteq \Nn(y)$. Therefore, the probability that $G$ is \emph{not} $(2,s)_j$-e.c.\ can be bounded from above by
\begin{align*}
n^{2+s} &\Big( 1 - (1-(1-p)^{|X|})(1-(1-p)^{|Y|}) (1-p)^{|\Nn[S]|} \Big)^{n-|\Nn[S]|-|\Nn(x,1)|-|\Nn(y,1)|} \\
&\le n^{s+2} \Big( 1 - (1+o(1)) (1-e^{-c^2})^2 \exp(-sp^2n) \Big)^{n-o(n)}  \\
&= \exp \Big( O( s \log n) - \Omega( \exp(-sp^2n) n ) \Big) \\
&= \exp \Big( O( \log^2 n) - \Omega( \sqrt{n} ) \Big) = o(1),
\end{align*}
provided that $sp^2n \le \frac 12 \log n$; that is, even for $s$ tending to infinity slowly enough the desired statement holds.

Now, let us consider the case $\log^3 n \le pn = o(\sqrt{n})$. Suppose that $(n \log n)^{\frac {1}{2j}} \ll pn \ll n^{1/j}$ for some $j = j(n) \in \N \setminus \{1\}$. (Note that for the range of $p$ under consideration there is at least one such $j$.) Let $S \subseteq V(G)$ ($s=|S|$), $x \in V(G) \setminus S$, and $y \in V(G) \setminus (S \cup \{x\})$. As before, Lemma~\ref{lem:gnp exp} and Lemma~\ref{lem:gnp exp2}(i) can be used (this time with $i=j-1$) to get that the size of the $(j-1)$-st neighbourhood of $S$ is $(1+o(1))s(pn)^{j-1}$, and that there are two disjoint sets $X$ and $Y$ that are also disjoint from the $(j-1)$st neighbourhood of $S$, each of cardinality $(1+o(1))(pn)^{j-1}$, such that $X \subseteq \Nn(x,j-1)$ and $Y \subseteq \Nn(y,j-1)$. This time, the probability that $G$ is \emph{not} $(2,s)_j$-e.c.\ can be bounded from above by
\begin{align*}
n^{2+s} &\Big( 1 - (1-(1-p)^{|X|})(1-(1-p)^{|Y|}) (1-p)^{|\Nn[S]|} \Big)^{n-|\Nn[S]|-|\Nn(x,j-1)|-|\Nn(y,j-1)|} \\
&\le n^{2s} \Big( 1 - (1+o(1)) (p(pn)^{j-1})^2 \Big)^{n-o(n)},
\end{align*}
since the condition $pn = o(n^{1/j})$ implies that $p(pn)^{j-1}=o(1)$. On the other hand, it follows from the condition $pn \gg (n \log n)^{1/(2j)}$ that $(p(pn)^{j-1})^2 n \gg \log n$. Hence, the probability is at most
$$
\exp \Big( O( s \log n) - \Omega( (p(pn)^{j-1})^2 n ) \Big) = o(1),
$$
provided that $s \log n \ll (p(pn)^{j-1})^2 n$; that is, even for $s$ tending to infinity slowly enough the desired statement holds. This finishes the proof.
\end{proof}

\subsection{Proof of Theorem~\ref{thm:sparse_graphs_nonconstant}}

We continue investigating sparse graphs but this time we would like to deal with $r-m$ tending to infinity with $n$ faster than before. In order to achieve this goal, we need to combine the two ideas used before and relax the $(2,s)_j$-e.c.\ property. For fixed positive integers $j$, $k$ and $l$, we say that a graph $G$ is $(1,l,k)_j$-\emph{existentially closed} (or $(1,l,k)_j$-\emph{e.c.}), if for any vertex $v \in V(G)$ and two disjoint subsets of $A, B \subseteq V(G) \setminus \{v \}$ with $|A|=l$ and $|B|=k$, there exists a vertex $z \in V(G) \setminus (\{v \} \cup A \cup B)$ at distance at most $j$ from $v$, at distance at most $j$ from \emph{some} vertex in $A$, and  at distance at least $j+1$ from every vertex in $B$. 
One final time, we have the following observation.

\begin{theorem}\label{thm:lkj_ec}
Let $j, r, m, s, l$ be positive integers such that $s \le r-m-l+1$, and let $G$ be any $(1,l,s)_j$-e.c.\ graph. Then
$$
\sigma(G,r,m) \ge s+1.
$$
\end{theorem}
\begin{proof}[Proof of Theorem~\ref{thm:sparse_graphs_nonconstant}]

Let $\eps=o(1)$ be a function tending to zero sufficiently slowly. Let $S \subseteq V(G) (s=|S|)$, $x \in V(G)\setminus S$, and $Y \subseteq V(G) \setminus (S \cup \{x\})$. In all the cases discussed below, we will set $s=(r-m)/(1-\eps)$, $|Y|=\eps s$, and our goal will be to show that there is some natural number $j$ such that a.a.s.\ $G$ is $(1,\eps s, s)_j$-e.c. By Theorem~\ref{thm:lkj_ec} this will then imply the result. We may assume that the properties mentioned in Lemmas~\ref{lem:gnp exp} and~\ref{lem:gnp exp2} are deterministically satisfied.

To prove the desired property for part (i), let us note that since $s=o(n/(pn)^2)$, it follows from Lemma~\ref{lem:gnp exp} (applied with $i=1$) that we may assume that the size of $\Nn[S]$ is $(1+o(1))spn$ (note also that by the same lemma $|\Nn(x,1)|=(1+o(1))pn$ and $|\Nn[Y]|=(1+o(1))\eps s p n$ a.a.s.). Also, by Lemma~\ref{lem:gnp exp2}(ii) (applied with $i=1$), there is a set $X$,  disjoint from the neighbourhood of $S$, such that  $|X| \geq (pn)/2$ and $X \subseteq \Nn(x)$. (Note that $s=o(pn)$, and thus the conditions of the Lemma are satisfied). Also, by applying Lemma~\ref{lem:gnp exp2}(ii) to each vertex $y \in Y$ separately (by adding $x$ and all other vertices of $Y$ that are already examined to the forbidden set $S$ in Lemma~\ref{lem:gnp exp2}(ii)), we obtain a subset $Z$ of size at least $(\eps s p n)/2$, where $Z \subseteq \Nn(Y)$, $Z$ being disjoint from $\Nn[S]$ and also disjoint from $\Nn(x)$. This time, the probability that $G$ is  \emph{not} $(1,\eps s,s)_2$-e.c.\ can be bounded from above by
\begin{align*}
n^{(1+\eps)s+1} &\Big( 1 - (1-(1-p)^{|X|})(1-(1-p)^{|Z|}) (1-p)^{|\Nn[S]|} \Big)^{n-|\Nn[S]|-|\Nn(x,1)|-|\Nn[Y]|} \\
&= \exp \Big( O( s \log n) - \Omega( p^2 n (\eps s p ^2 n) n ) \Big) \\
&= \exp \Big( O( s \log n) - \Omega( \eps s p^4 n^3) \Big) \\
&\leq \exp \Big( s  \Big(  O(\log n) - \Omega(\eps \omega \log n) \Big) \Big) = o(1),
\end{align*}
where the first equality follows from the observation that $p|\Nn[S]| =O(p^2 s n)=o(1)$ (which is equivalent to $s=o(n/(pn)^2))$, and the last line follows from the condition $pn \geq (\omega n \log n)^{1/4}.$ Thus, a.a.s.\ $G$ is $(1,\eps s, s)_2$-e.c.\ 

To show the desired property for part (ii), using the notation from before, it follows from Lemma~\ref{lem:gnp exp} (applied with $i=2$) that  $|\Nn[S,2]|=(1+o(1))s(pn)^2$. (Note that since $s=o(pn)=o(\omega n \log n)^{1/4}$ one can apply the lemma with $i=2$.) Also, by Lemma~\ref{lem:gnp exp2}(ii) applied with $i=2$ (note that the condition $2\eta < 1-\eta$ from the lemma is satisfied), there are two disjoint sets $X$ and $Z$ that are also disjoint from $\Nn[S,2]$, $X$ of cardinality at least $(pn)^2/2$, $Z$ of cardinality at least $(\eps s) (p n)^2/2$, such that $X \subseteq \Nn(x,2)$, $Z\subseteq \Nn[Y,2]$, and $Z$ disjoint from $\Nn(x,2)$. Since $p|\Nn[S,2]|=O(p^3 s n^2)=o(1)$ by the conditions $s=o(pn/\omega \log n)$ and $pn \leq (\omega n \log n)^{1/4}$,  the probability that $G$ is not $(1,\eps s,s)_3$-e.c.\ is 
\begin{align*}
n^{(1+\eps)s+1} &\Big( 1 - (1-(1-p)^{|X|})(1-(1-p)^{|Z|}) (1-p)^{|\Nn[S,2]|} \Big)^{n-|\Nn[S,2]|-|\Nn(x,2)|-|\Nn[Y,2]|} \\
&= \exp \Big( O( s \log n) - \Omega( p^3 n^2 ((\eps s) p ^3 n^2) n ) \Big) \\
&= \exp \Big( s \Big( O(\log n) - \Omega( \eps p^6 n^5) \Big)  \Big) \\
&= o(1),
\end{align*}
where the last line follows as before by plugging in the lower bound for $pn$. Thus, a.a.s.\ $G$ is $(1,\eps s, s)_3$-e.c.\ 

Finally, in order to prove Part (iii) let us take the unique constant $j \geq 4$ such that $(\omega n \log n)^{\frac{1}{2j}} \leq pn < (\omega n \log n)^{\frac{1}{2(j-1)}}$ (recall that $pn < (\omega n \log n)^{1/6}$ and $pn = n^{\eta+o(1)}$ for some $\eta \in (0,1/2)$). Then, analogously to Part (ii), one can apply Lemma~\ref{lem:gnp exp} and Lemma~\ref{lem:gnp exp2}(ii), but this time with $i=j-1$. Since $p|\Nn[S,j-1]| \leq ps(pn)^{j-1}=o(1)$, by the same calculations as in Part (ii), the probability that $G$ is  not $(1,\eps s,s)_{j}$-e.c.\ is $o(1)$. The proof is complete.
\end{proof}


\begin{thebibliography}{99}

\bibitem{anw} M.H.~Albert, R.J.~Nowakowski, and D.~Wolfe, \emph{Lessons in Play}, A K Peters, Ltd., 2007.

\bibitem{al} B.\ Alspach, Sweeping and searching in graphs: a brief survey, \emph{Matematiche} \textbf{59} (2006) 5--37.

\bibitem{AB1} A.~Bonato and R.J.\ Nowakowski, \emph{The Game of Cops and Robbers on Graphs}, American Mathematical Society, Providence, Rhode Island, 2011.


\bibitem{bkl} B.\ Bollob\'as, G.\ Kun, I.\ Leader, Cops and robbers in a random graph, preprint.

\bibitem{bpw} A.\ Bonato, P.\ Pra\l{}at, C.\ Wang, Network Security in Models of Complex Networks, \emph{Internet Mathematics} \textbf{4} (2009), 419--436.

\bibitem{BCPWZ} J.V.~Butterfield, D.W.~Cranston, G.~Puleo, D.B.~West, and R.~Zamani, Revolutionaries and Spies: Spy-good and spy-bad graphs,  \emph{Theoretical Computer Science} \textbf{463} (2012), 35--53.

\bibitem{CSW} D.W.~Cranston, C.D.~Smyth, and D.B.~West, Revolutionaries and spies in trees and unicyclic graphs, \emph{Journal of Combinatorics} \textbf{3} (2012), 195--206.

\bibitem{fomin} F.V.\ Fomin and D.\ Thilikos, An annotated bibliography on guaranteed graph searching, \emph{Theoretical Computer Science} \textbf{399} (2008) 236--245.

\bibitem{gena} G.\ Hahn, Cops, robbers and graphs, \emph{Tatra Mountain Mathematical Publications} \textbf{36} (2007) 163--176.

\bibitem{HS} D.~Howard and C.D.~Smyth, Revolutionaries and spies on grid-like graphs, preprint.

\bibitem{JLR} S.\ Janson, T.\ {\L}uczak, A.\ Ruci\'nski, \emph{Random Graphs}, Wiley, New York, 2000.

\bibitem{lp2} T.~\L{}uczak and P.~Pra\l{}at, Chasing robbers on random graphs: zigzag theorem, \emph{Random Structures and Algorithms} \textbf{37} (2010), 516--524. 

\bibitem{nw} R.\ Nowakowski, P.\ Winkler, Vertex to vertex pursuit in a graph, \emph{Discrete Mathematics} \textbf{43} (1983) 230--239.

\bibitem{p} P.\ Pra\l{}at, When does a random graph have constant cop number?, \emph{Australasian Journal of Combinatorics} \textbf{46} (2010), 285--296.

\bibitem{pw} P.\ Pra\l{}at and N.\ Wormald, Meyniel's conjecture holds for random graphs, preprint.

\bibitem{pw2} P.\ Pra\l{}at and N.\ Wormald, Meyniel's conjecture holds for random $d$-regular graphs, preprint.

\bibitem{q} A. Quilliot, Jeux et pointes fixes sur les graphes, Ph.D. Dissertation, Universit\'{e} de Paris VI, 1978.

\end{thebibliography}
\end{document}